\documentclass[11pt]{article}
\usepackage[T1]{fontenc}
\usepackage{lmodern,amsmath,amssymb}
\usepackage[a4paper,left=.92in,right=.92in,top=.84in,bottom=.84in]{geometry}
\usepackage[hidelinks]{hyperref}
\hypersetup{pdftitle={Special Smith forms of the meander and chromatic-join matrices},pdfauthor={Iaroslav Argunov}}
\newtheorem{theorem}{Theorem}
\newtheorem{lemma}[theorem]{Lemma}
\newtheorem{corollary}[theorem]{Corollary}
\newenvironment{proof}[1][Proof]{\par\noindent\textit{#1.}\ }{\unskip\nobreak\hfill\(\square\)\par}
\newcommand{\ZZ}{\mathbb Z}
\newcommand{\F}{\mathbb F}

\newcommand{\nc}{\mathcal P_n}
\newcommand{\rk}{\operatorname{rank}}

\begin{document}
\begin{center}
{\LARGE\bfseries Special Smith forms of the meander\\[.16em] and chromatic-join matrices\par}
\vspace{.7em}{\large Iaroslav Argunov\par}
\smallskip Independent Researcher, Oakland, USA\par
\vspace{.2em}Revised November 24, 2026
\end{center}
\begin{abstract}
We construct a special Smith form over \(\mathbb Z[X,Y]\) for the
two-variable chromatic-join matrix of noncrossing partitions, settling
the first conjecture of Miller and Stanton. Their second conjecture,
for the meander Gram matrix over \(\mathbb Z[x]\), follows by their
equal-parameter specialization. A finite B\'ezout identity combines
integral multiples of two-coloured Jones--Wenzl projectors and produces
a unitriangular operation on rows of each Dyck height. In the resulting
factorization \(LJ_n=DV\), the entries of the polynomial quotient \(V\)
are determined by two block counts, the row height and a prefix of arcs in distinct cycles.
A reflection
count accounts for the full determinant degree and makes the quotient
unimodular. The connected matrix has integer Smith factors determined
by the maximum heights of the paths.
\end{abstract}
\section{Matchings and partitions}

Write \(\nc\) for the noncrossing perfect matchings of \(2n\) labelled
points in boundary order, \(n\ge1\). The union of \(P,Q\in\nc\)
consists of \(c(P,Q)\) alternating cycles, with a common edge counted
as a cycle. The meander Gram matrix is
\[
                 G_n(x)_{P,Q}=x^{c(P,Q)} .
\]
The row and column orders agree. The heights of \(P\) are the numbers
of arcs intersecting the gaps in the linear order of its endpoints.
They give a Dyck path of length \(2n\); let \(h(P)\) be its maximum.
For \(d\ge2\) define the monic polynomial
\[
 \Psi_d(x)=\prod_{\substack{1\le j<d\\(j,d)=1}}
                  \bigl(x-2\cos(\pi j/d)\bigr),\qquad
 s_k(x)=\prod_{d=2}^{k+1}\Psi_d(x).
\]
A special Smith form is a diagonal form with successive entries
dividing one another, obtained by left and right transformations of
determinant \(1\). The ring \(\ZZ[x]\) is not a principal ideal domain.

For the cyclotomic polynomial \(\Phi_d\), set \(x=w+w^{-1}\) and \([1]=1,[2]=x,[j+1]=x[j]-[j-1]\). Then
\begin{equation}
 w^{\varphi(d)}\Psi_d(w+w^{-1})=\Phi_d(w^2),\qquad
 [j]=\prod_{\substack{d\mid j\\d>1}}\Psi_d . \label{cyc}
\end{equation}
The first identity follows by matching roots and leading coefficients.
Its Laurent expression after division by \(w^{\varphi(d)}\) is symmetric
in \(w,w^{-1}\) with integer coefficients, hence lies in \(\ZZ[w+w^{-1}]\).
The second identity follows by factoring
\([j](w+w^{-1})=w^{1-j}(1+w^2+\cdots+w^{2j-2})\).

\section{The two-variable chromatic-join form}
\label{sec:two-colour}

Write \(S=\ZZ[X,Y]\) and \(c=(0\,1\,\dots\,n-1)\). For a noncrossing
partition \(\pi\) of \(\{0,\dots,n-1\}\), denote also by \(\pi\) the
permutation obtained by reading each block in increasing cyclic order.
Its Kreweras complement is the partition of cycles of \(\pi^{-1}c\),
denoted by \(\pi'\). Joins in the following matrix are taken in the
lattice of \emph{all} set partitions, and vertical bars count blocks:
\[
 J_n(X,Y)_{\pi,\rho}=X^{|\pi\vee\rho|}Y^{|\pi'\vee\rho'|}.
\]
Put \(B_2(z)=z\), and for \(d>2\) define \(B_d\in\ZZ[z]\) by
\(\Psi_d(w)=B_d(w^2)\). The latter equality follows by pairing the
opposite roots of the monic even polynomial \(\Psi_d\). These are
the factors \(f_{d-1}\) of \cite[\S6.2.2]{ms}. For \(k\ge1\) put
\[
                     b_k(X,Y)=\prod_{d=2}^{k+1} B_d(XY).
\]

\begin{theorem}\label{thm:two-colour}
Over \(\ZZ[X,Y]\), the special Smith form of \(J_n(X,Y)\) has one
diagonal entry \(b_{h(P)}(X,Y)\) per noncrossing matching \(P\),
in order of nondecreasing maximum height. This proves
Conjecture~1 of Miller and Stanton.
\end{theorem}

First identify the entries with a spherical evaluation. Mark the
open boundary interval from \(2i\) to \(2i+1\) black and the interval
from \(2i+1\) to \(2i+2\) white, with endpoints taken modulo \(2n\).
The matching attached to \(\pi\) joins \(2i+1\) to \(2\pi(i)\).
It is the boundary of a thin neighbourhood of the disjoint polygons
of the blocks, cut at the black intervals; conversely, the black
regions recover these blocks. Following a white interval \(i\)
through its endpoint \(2i+2\) reaches the odd endpoint indexed
by \(\pi^{-1}(i+1)\). Thus the white regions give precisely \(\pi'\).
After two disks are glued along their labelled boundaries, every
complementary region meets a boundary interval, and the equivalence
relation among intervals of either colour is generated by the
relations from the two disks. Assigning weight \(X\) to every black
region and \(Y\) to every white region on the sphere gives
\(J_n(X,Y)_{\pi,\rho}\). In particular every entry is divisible by \(XY\).

We use Hazi's two-coloured Temperley--Lieb category
\cite[Theorem~A, equation~(7), Definition~3.3 and Lemma~3.5]{hazi}.
Assign the parameters \(-X,-Y\) so that the colour inside a removed disk receives its named weight. In Hazi's convention a blue disk inside a red region has weight \(-[2]_s\). In these conventions deletion of an
innermost circle has weight \(X\) or \(Y\) according to its interior
colour. The planar evaluation omits its exterior region: to obtain the
spherical weight, multiply by the colour of this last region. This
also shows directly that the spherical evaluation is independent of
which complementary region is chosen as exterior.

Write \(\Theta_{d,u}\) for Hazi's primitive factor of index \(d\)
in colour \(u\). For \(d=2\) it is the corresponding parameter
\(-X\) or \(-Y\); for \(d>2\) it is \(B_d(XY)\) for both colours.
The two-colour quantum integer is
\([j]_u=\prod_{d\mid j,\,d>1}\Theta_{d,u}\).
Write \(v_P\) for the row indexed by \(P\). Its area is the sum of
the Dyck heights over all gaps. We use two facts about the generic
two-colour projector: its coefficient on the identity diagram is one,
and multiplication by an adjacent turnback annihilates it.
Hazi's denominator and trace formulas will specify the multiple
that can be inserted integrally. The argument below works in the
polynomial ring even at specializations where both variables are
nonunits. Earlier determinant formulas are due to Dahab \cite{dahab}
and Di Francesco, Golinelli and Guitter \cite{dfgg}.
Chen and Przytycki related chromatic joins to ordinary Temperley--Lieb
pairings by a diagonal rescaling \cite{chen}. For categories with coloured
regions, Howat, Laugwitz and Ray obtained Gram determinant and
nondegeneracy formulas over a field \cite[Corollary~4.21]{hlr}.
Here we seek a diagonal equivalence over the polynomial ring itself.
The projector existence criterion over a ring is discussed in
\cite{webster,hazi}.

\section{Integral elimination by a finite B\'ezout identity}
\label{sec:integral-elimination}

There is a direct construction of the row operation for \(J_n\).
Its coefficients use a single univariate B\'ezout identity for each height.
Write \(z=XY\), \(N=k+1\), and retain \(B_2(z)=z\) and the
factors \(B_d(z)\) defined above. For \(2\le r\le N\) put
\[
 t_{N,r}(z)=\prod_{\substack{2\le d\le N\\d\nmid r}}B_d(z),
 \qquad
 I_N=\{r:2\le r\le N,\ 2r>N\}.
\]

\begin{lemma}
There exist \(a_{N,r}(z)\in\mathbb Z[z]\) such that
\(\sum_{r\in I_N}a_{N,r}(z)t_{N,r}(z)=1\).
\end{lemma}
\begin{proof}
For incomparable indices \(d,e\) under divisibility the ideal
\((B_d,B_e)\) is the unit ideal of \(\mathbb Z[z]\).
Indeed specialize \(x_s=z,x_t=1\) in Hazi's
two-colour B\'ezout identity for the primitive factors of colour \(s\)
\cite[Lemma 3.5]{hazi}. This also covers an index equal to two.
Multiply one such identity over every incomparable pair from
\(\{2,\ldots,N\}\). In each expanded summand the indices whose
factors were never selected form a chain. Choose \(r\in I_N\) divisible
by the largest member of the chain, or any \(r\in I_N\) if the chain
is empty. Such a multiple exists by choosing the largest multiple
up to \(N\). Each \(d\nmid r\) has therefore been selected at least
once, so the summand is a multiple of \(t_{N,r}\). Group the summands.
The same reasoning applies to the empty product of identities.
\end{proof}

\begin{lemma}
For a row of height \(k\) and each \(r\in I_N\), one can add an
integral polynomial combination of rows of strictly smaller area
to \(t_{N,r}(XY)v_P\) to obtain an element of
\(b_k(X,Y)\mathbb Z[X,Y]^{|\mathcal P_n|}\).
\end{lemma}
\begin{proof}
Work first in the generic two-colour diagram calculus. For either
orientation \(u\), a common denominator of the coefficients of
the Jones--Wenzl projector \(f^u_{r-1}\) is
\[
 g_{r-1,u}=\prod_{\substack{2\le d<r\\d\nmid r}}\Theta_{d,u}.
\]
This is the denominator statement in Hazi's proof of Theorem A
\cite[\S4, Lemma 4.2 and the proof of Theorem A]{hazi}.
Under \(x_s=-X, x_t=-Y\) the factors with \(d>2\) give
\(B_d(XY)\). An index \(d=2\) contributes either \(-X\) or \(-Y\),
both of which divide \(B_2(XY)\). Hence \(g_{r-1,u}\mid t_{N,r}(XY)\).
The box \(t_{N,r}(XY)f^u_{r-1}\) is an integral polynomial
linear combination of diagrams, with identity coefficient \(t_{N,r}(XY)\).

Insert this box on \(r-1\) consecutive arcs crossing a maximal gap in
\(P\). Every nonidentity term reduces the number crossing that gap.
At any earlier gap, only the selected endpoints already encountered
can contribute a crossing of a redrawn arc, and each contributed one
previously. Use the reversed order at later gaps. Thus the remaining
terms have smaller area. Each selected half-arc leads to a distinct
marked boundary endpoint, so expanding the box creates no closed
circle on the row side alone.

Fix the leftmost gap of maximum height in each row \(P\). At that
gap order the crossing arcs from outermost to innermost. Use the
first \(m=r-1\) for the insertion, so that the choices for different
\(r\) are nested. Write \(K_{P,Q,r}\) for the value of the entry
after inserting \(t_{N,r}(XY)f^u_m\) and closing with \(Q\).
Put \(z=XY\) and
\[
                 E_r(z)=\prod_{\substack{d\mid r\\ d>2}} B_d(z).
\]

The complement of the box in the sphere is a disk. It contains
circles and a matching of the \(2m\) ports. Delete the circles starting
with any innermost disk on the side away from the box. Such a disk
contains no port or arc ending at a port. The deletion contributes
the colour weight of its interior face. The resulting monomial \(C\)
is common to all terms of the projector and to the identity box:
each term has the same prescribed boundary shading. After all these
deletions, an arc returning to the same row of ports has an innermost
adjacent turnback on that row, so the value with the projector is
zero. This uses the turnback relation over the fraction field;
the cleared box already has coefficients in the polynomial ring.

If every complementary arc passes from one row of ports to the other,
noncrossing forces the unique order-preserving matching. Choose
the exterior point next to the left edge of the box and draw this
matching as a right closure. Let \(w\in\{X,Y\}\) be the weight of the
left face, and \(w'\) the opposite weight, so that \(ww'=z\).
For the left boundary colour \(u\), Hazi's disk rule gives
\([2]_u=-w'\). Closing the rightmost strand keeps the left boundary
colour fixed. Iterating the partial trace identity
\(\operatorname{pTr}(f_j^u)=-[j+1]_u f_{j-1}^u/[j]_u\)
in the fraction field gives the planar scalar \((-1)^m[r]_u\)
\cite[equation (7)]{hazi}. Multiplying by \(w\) counts the outer
face once. Section~\ref{sec:spherical-check} derives the same scalar
from the ordinary projector. The spherical values of the boxes alone are therefore
\[
 W_f=
 \begin{cases}
       z E_r(z),&r\ \text{even},\\
       w E_r(z),&r\ \text{odd},
 \end{cases}
 \qquad
 W_{\rm id}=w^{\lfloor m/2\rfloor+1}(w')^{\lceil m/2\rceil}.
\]
Indeed \([r]_u=E_r(z)\) for odd \(r\) and
\([r]_u=-w' E_r(z)\) for even \(r\). For the identity, the \(m\)
nested circles give alternating face colours starting with \(w\).
Choosing the exterior point on the other side gives the same two
spherical values: the side colours coincide for even \(m\), and both
expressions are symmetric in them for odd \(m\). The boundary ports
of the box stay fixed throughout. With the outside of the stack on
the left, \(w=Y\) for the chosen prefix: the outer face of the
matching meets the white interval from \(2n-1\) to \(0\).

Write \(a=|\pi\vee\rho|\) and \(a'=|\pi'\vee\rho'|\). Restoring
the identity box recovers the original pair, hence
\(C W_{\rm id}=J_n(P,Q)=X^aY^{a'}\). Combining the two values with
the definitions of \(t_{N,r}\) and \(b_k\) gives the exact quotient
\begin{equation}
\frac{K_{P,Q,r}}{b_k(X,Y)}
=
\begin{cases}
 X^{a-\lfloor r/2\rfloor}Y^{a'-\lfloor r/2\rfloor},
   &\text{the selected arcs lie in distinct cycles of }P\cup Q,\\
 0,&\text{otherwise}.
\end{cases} \label{eq:through-quotient}
\end{equation}
To see the criterion, cut the \(m\) selected arcs in the identity
diagram. When they belong to distinct cycles, the two ends of each
one are connected by the remainder of its own cycle. Conversely, the
unique through matching reconnects precisely those pairs, so
restoring the identity puts each selected arc on its own cycle.
In the through case the displayed quotient equals \(C\) for even
\(r\) and \(wC\) for odd \(r\); its exponents are nonnegative.
This proves the required integral divisibility for every pair.

\end{proof}

Multiply the rows of the second lemma by the coefficients of the
first lemma and add. The coefficient of \(v_P\) is exactly one.
Doing this for every \(P\), in increasing order of area, constructs
a lower unitriangular \(L\) with \(LJ_n=DV\), where \(V\) has
integral polynomial entries. The degree count in Section~\ref{sec:degree}
gives \(\det V=1\). In the Smith diagonalization \(LJ_n R=D\) the
right multiplier is \(R=V^{-1}\). Simultaneous signed permutations
sort the heights.
For \(k=3\), the scalar identity in this construction is simply
\[
             1=(z-1)^2-z(z-2),\qquad I_4=\{3,4\}.
\]
Finally, at \(X=Y=x\) the spherical region count is the loop count
plus one. Thus \(J_n(x,x)=xG_n(x)\), while
\(b_k(x,x)=x\,s_k(x)\). Cancelling \(x\) in the matrix identity
gives the one-variable special Smith form, preserving both
determinants. Explicitly, its entries are \(s_{h(P)}(x)\) in
nondecreasing order of height. This is Conjecture~2. The implication
from the two-variable form was already recorded by Miller and
Stanton \cite[arXiv version, equation (64)]{ms}.

\begin{corollary}\label{cor:entries}
Fix the row order, the outermost prefixes and the B\'ezout weights above.
The entries of the quotient \(V\) in \(LJ_n=DV\) are determined
by four integers. For \(P,Q\) write \(k=h(P)\), and let \(q\) be
the length of the longest initial segment of the \(k\) chosen
crossing arcs that belong to pairwise distinct cycles of \(P\cup Q\).
Let \(a,a'\) be the two block counts in \(J_n(P,Q)\) and set
\(\ell=\lfloor(q+1)/2\rfloor\). Then
\[
 V_{P,Q}=X^{a-\ell}Y^{a'-\ell}\,
 \sum_{\substack{r\in I_{k+1}\\r\le q+1}}
   a_{k+1,r}(XY)(XY)^{\,\ell-\lfloor r/2\rfloor}.
\]
\end{corollary}
\begin{proof}
In the weighted sum of \eqref{eq:through-quotient} the surviving
prefixes have length at most \(q\). Applying the identity-box face
count to the first \(q\) arcs gives \(a,a'\ge\ell\), whether or not
\(q+1\) belongs to \(I_{k+1}\). Factoring out the displayed
monomial gives the result.
\end{proof}
Thus at order \(n\) one can precompute at most \(n(n+1)/2\) univariate
polynomials, one for each pair \(1\le q\le k\le n\).
For a particular pair of matchings, a traversal of its alternating
cycles determines \(q\). If the total number of cycles is less than
\(\lfloor(k+1)/2\rfloor\), the sum is empty. In particular a connected
pair gives a zero entry of \(V\) whenever the row height is at least
three. These statements concern the change of basis.

\section{The degree of the determinant}
\label{sec:degree}
One can finish the unimodularity argument without assuming a determinant
formula. We give the calculation for both rings. Write
\(C_n=\binom{2n}{n}/(n+1)\) and \(c_j=\binom{2n}{n+j}\), extended by
zero outside its range.

\begin{lemma}
\(\displaystyle\sum_{P\in\nc}\deg s_{h(P)}=nC_n.\)
\end{lemma}
\begin{proof}
The reflection formula for walks between the absorbing walls at
\(-1\) and \(d-1\) gives
\(\sum_{j\in\mathbb Z}(c_{jd}-c_{jd+1})\) for the number of
Dyck paths of length \(2n\) with maximum at most \(d-2\).
For completeness, the image sum at an intermediate endpoint \(i\)
uses the unrestricted counts from each \(2jd\) minus those from
\(-2+2jd\). It vanishes at both walls by pairing reflected images,
satisfies the nearest-neighbour recurrence, and at time zero is the
point mass at the origin in the allowed strip. Evaluating it at
time \(2n\) and endpoint zero gives the stated formula.
Pair the terms with indices \(j\) and \(-j\). The number reaching
height \(d-1\) is
\(\sum_{m\ge1,\,d\mid m}(c_{m-1}-2c_m+c_{m+1})\).
Since \(\deg\Psi_d=\varphi(d)\), summing these counts and using
\(\sum_{d\mid m,\,d>1}\varphi(d)=m-1\) gives
\[
 \sum_{P\in\nc}\deg s_{h(P)}
 =\sum_{m=2}^{n+1}(m-1)(c_{m-1}-2c_m+c_{m+1})
 =c_1=nC_n.
\]
\end{proof}

The diagonal of \(G_n\) has monomials \(x^n\); every off-diagonal
entry has smaller degree. Thus its determinant is monic of degree
\(nC_n\). Once \(L G_n=D V\) is integral with \(L\) unitriangular,
the degree calculation forces \(\det V=1\).

Likewise each diagonal entry of \(J_n\) has total degree \(n+1\):
the sphere has one more region than loops. Every off-diagonal entry
has total degree at most \(n\), so the unique determinant term of
maximal degree comes from the identity permutation. Each \(b_h\)
is a monic polynomial in the product \(XY\), so setting \(X=Y=x\)
preserves its total degree. From \(b_h(x,x)=x\,s_h(x)\) this degree is
\(1+\deg s_h\).
Consequently the integral quotient determinant in \(LJ_n=DV\) is
constant. Setting \(X=Y=x\) compares the monic leading terms, and
the constant is \(+1\). This recovers the product of diagonal factors
in Dahab's determinant evaluation.

\section{The connected matrix}

Set \(H_n(P,Q)=\mathbf1_{\{c(P,Q)=1\}}\). Each entry of \(G_n(x)\)
is divisible by \(x\), and \(H_n=(G_n(x)/x)|_{x=0}\).
\begin{corollary}\label{integer}
The integer Smith factors of \(H_n\), up to the usual positive signs,
are
\[
       \operatorname{lcm}\bigl(1,\ldots,\lfloor(h(P)+1)/2\rfloor\bigr),
           \qquad P\in\nc .
\]
In particular the rank over \(\F_p\) is the number of length-\(2n\)
Dyck paths bounded above by \(2p-2\).
\end{corollary}
\begin{proof}
Divide the polynomial matrix identity through by \(x=\Psi_2(x)\)
and then set \(x=0\); the unimodular transformations specialize.
With \(w=i\), \eqref{cyc} gives
\(\left|\Psi_d(0)\right|=|\Phi_d(-1)|\).
For \(d>2\) this equals \(p\) when \(d=2p^a\) for a prime \(p\)
and \(a\ge1\), and equals \(1\) otherwise. This follows from
\(\Phi_{2j}(-1)=\Phi_j(1)\) for odd \(j>1\), the elementary
cyclotomic values at \(1\), and the power-of-two case.
Multiplication gives the displayed least common multiple.
It is divisible by \(p\) exactly when \(h(P)\ge2p-1\).
\end{proof}

For \(n=5\) the 42 paths give integer factors \(1^{16},2^{25},6^1\);
in particular the connected matrix has rank 41 in characteristic
three. The full polynomial equivalence also specializes at every
integer loop weight, including the singular values. In the notation
of Miller and Stanton, \(J_n(x,x)=xG_n(x)\), so their equal-parameter
chromatic-join specialization follows as well.

\section{The fixed-block-count grove pairing}
\label{sec:groves}

Let \(1\le b\le n\), and index the rows of \(H_{n,b}\) by the
noncrossing partitions \(\pi\) with \(b\) blocks and its columns by
those \(\rho\) with \(n+1-b\) blocks. Put \(H_{n,b}(\pi,\rho)=1\)
if \(\pi\vee\rho\) is the one-block partition and \(0\) otherwise.
These square matrices are the complementary blocks of the grove
pairing of Kenyon and Wilson
(\href{https://arxiv.org/abs/math/0608422}{\emph{Boundary partitions
in trees and dimers}}, \S2.2, equation (3)).

\begin{corollary}\label{cor:grove-block}
The positive Smith factors of \(H_{n,b}\) are
\[
 \operatorname{lcm}\bigl(1,\ldots,\lfloor(h(P)+1)/2\rfloor\bigr)
\]
over the Dyck paths \(P\) of length \(2n\) with exactly \(b\) upsteps
in odd-numbered positions, counting the first step as position one.
In particular \(\rk_{\F_2}H_{n,b}=\binom{n-1}{b-1}\).
\end{corollary}

\begin{proof}
Write \(a=|\pi\vee\rho|\) and \(a'=|\pi'\vee\rho'|\).
The plane bipartite multigraph with one vertex for each block of
\(\pi\) and \(\rho\), and one edge per label joining its two blocks,
has \(a\) components and \(a'\) faces. To see the latter, draw the two
families of stars in opposite hemispheres. Their complementary regions
join the white boundary intervals exactly as in
Section~\ref{sec:two-colour}. Euler's formula gives
\[
             a-a'=|\pi|+|\rho|-n-1 .
\]
Grade \(\ZZ[X,Y]\) by \(\deg X=1,\deg Y=-1\). The entries of
\(J_n/(XY)\) have the degrees on the right, and all diagonal factors
\(b_{h(P)}/(XY)\) have degree zero. In the row relation of
Theorem~\ref{thm:two-colour}, retain from the coefficient of each
smaller row \(\tau\) only its homogeneous part of degree
\(|\pi|-|\tau|\). Coordinatewise this extracts the homogeneous
part of the degree of the original entry. The remainder is still
divisible by the same homogeneous diagonal factor. We thus obtain
\[
 L\,(J_n/(XY))=D V,\qquad D_{P,P}=b_{h(P)}(X,Y)/(XY),
\]
where \(L\) is unitriangular by area and
\(\deg V_{\pi,\rho}=|\pi|+|\rho|-n-1\) whenever the entry is nonzero.
The determinant argument gives \(\det V=1\) for these choices as well.

At \(X=Y=0\) a coordinate of \(J_n/(XY)\) survives precisely when
\(a=a'=1\). The displayed degree identity identifies its complementary
blocks with \(H_{n,b}\). Only degree-zero entries of \(L,V\) can
survive, so \(L(0,0)\) preserves block count and \(V(0,0)\) pairs
complementary counts. Each complementary square block of \(V(0,0)\)
has determinant \(\pm1\), since their product up to sign is a unit.
Each block of \(L(0,0)\) is unitriangular. Finally
\(\lvert D_{P,P}(0,0)\rvert\) is the least common multiple in the
statement, by \(B_d(0)=\Psi_d(0)\) for \(d>2\) and
Corollary~\ref{integer}. In the fattening convention of
Section~\ref{sec:two-colour}, the step at endpoint \(2i\), numbered
\(2i+1\), opens an arc exactly when \(i\) is the minimum of its
block. This identifies the required paths.
\end{proof}

More generally, put \(u=p^e\) for a prime \(p\) and \(e\ge1\).
Let \(T_u(z)\) be the \((u-1)\)-by-\((u-1)\) tridiagonal matrix
with diagonal \(1+z\), superdiagonal \(z\) and subdiagonal \(1\).
The number of Smith factors of \(H_{n,b}\) having \(p\)-adic
valuation less than \(e\) is
\[
              [z^{b-1}]\bigl(T_u(z)^{n-1}\bigr)_{1,1}.
\]
Indeed the height cutoff from Corollary~\ref{cor:grove-block} is
\(2u-2\). Delete the first and last steps, then read the remaining
steps in pairs between heights \(1,3,\ldots,2u-3\). Marking the
second step of each pair by \(z\) if it rises gives precisely
the entries of \(T_u\); the deleted first step supplies the one
remaining marked rise. In particular \(T_2=(1+z)\), giving the
binary rank in the corollary.

For \(n=5,b=3\), the integral diagonal has six entries \(1\),
thirteen entries \(2\), and one entry \(6\). A direct Smith reduction
of the twenty-by-twenty matrix of joins gives the same list. The
corresponding multiplicities for \(n=8,b=4\) are
\(1^{35},2^{322},6^{126},12^7\).

\section{The spherical normalization of the projector}
\label{sec:spherical-check}

We give a second derivation of the scalar used in
\eqref{eq:through-quotient}. It also distinguishes a change of
exterior region from a change of the prescribed colour of a port.
Use \(W\) for the weight of the region at the left edge of the box
and \(U\) for the other weight. Work in the generic field with
\(t^2=WU\), and write \(q_0=0,\ q_1=1,\
q_{j+1}=tq_j-q_{j-1}\). Write \(I\) for the matching by straight
strands and also for the complementary matching closing each strand
back to itself.

\begin{lemma}
The spherical evaluation of the two-colour projector on \(m\ge1\)
strands against \(I\), with the exterior face counted, is
\[
 S_m(W,U)=
 \begin{cases}
 Wq_{m+1}(t),&m\ \mathrm{even},\\
 tq_{m+1}(t),&m\ \mathrm{odd}.
 \end{cases}
\]
Every other exterior matching evaluates to zero.
\end{lemma}
\begin{proof}
Call the shade of weight \(W\) black. For disk matchings \(P,Q\),
write \(r(P),r(Q)\) for the numbers of their black regions before
gluing. Suppose their spherical gluing has \(\ell\) circles, \(a\)
black regions and \(b\) white regions. Then \(a+b=\ell+1\) and
\[
 a-b=r(P)+r(Q)-m-1.
\]
Indeed, place one vertex in each black region of each disk and join
the two incident vertices across each black boundary interval. The
plane graph has \(r(P)+r(Q)\) vertices, \(m\) edges, \(a\) components
and \(b\) faces. This is the region count behind the diagonal
rescaling in \cite{chen}.

Set \(\eta=W/t,\ h(P)=\eta^{r(P)-r(I)}\), and
\(\alpha=t\eta^{2r(I)-m-1}\). The two face counts give
\[
 W^aU^b=\alpha h(P)h(Q)t^{\ell}.
\]
The parallel strands in \(I\) separate alternating regions starting
with black, hence \(r(I)=1+\lfloor m/2\rfloor\). Thus
\(\alpha=W\) for even \(m\) and \(\alpha=t\) for odd \(m\).

Write the ordinary projector as \(\sum_P c_P P\), with \(c_I=1\).
Every exterior matching \(Q\ne I\) has an innermost adjacent
same-row pair. The turnback equations kill its pairing with the
ordinary projector. By the rescaling, the shaded row with coefficients
\(c_P/h(P)\) has the same vanishing property as the shaded projector.
The ordinary Gram matrix is invertible over the generic field:
the diagonal has degree \(m\) and all other entries have smaller
degree in \(t\). Its rescaling remains invertible. Omitting the column
indexed by \(I\) therefore leaves a one-dimensional left annihilator.
Both rows in it have coefficient one at \(I\), so they agree.
Their evaluation against \(I\) is \(\alpha\) times the ordinary trace.

For completeness, write \(f_m\) for the ordinary projector,
\(P=f_{m-1}\otimes1\), and \(e=e_{m-1}\). In the positive loop
convention \(e^2=te\), Wenzl's recurrence gives
\(f_m=P-(q_{m-1}/q_m)PeP\) for \(m\ge2\)
\cite[\S2, with the opposite loop sign]{morrison}.
The full closures of \(P\) and \(PeP\) have respective traces
\(t\operatorname{tr}(f_{m-1})\) and \(\operatorname{tr}(f_{m-1})\),
by cyclicity, \(P^2=P\), and the cap and cup at the last strand.
Starting from \(\operatorname{tr}(f_1)=t\) gives
\(\operatorname{tr}(f_m)=q_{m+1}(t)\).
\end{proof}

The right closure joins each upper port to the lower port immediately
below it by nested arcs around the right edge of the rectangle. The
face at the left edge retains weight \(W\). For instance
\(S_2(W,U)=W(WU-1)\), so \(S_2(2,3)=10\) and \(S_2(3,2)=15\),
whereas \(S_1(W,U)=WU\) is symmetric. These scalars lie in
\(\mathbb Z[W,U]\) by parity of the recurrence. The use of a generic
field here is confined to the trace calculation. The integral
denominator is established in Section~\ref{sec:integral-elimination}.

Let \(R\) be the fully nested matching and \(J_m^{\widehat R}\)
the matrix obtained by deleting its row and column. In the cyclic
port order the straight strands give \(R\), and the gap before the
first port is white. If \(v\) is the coefficient row of the projector,
then \(v_R=1\) and \(vJ_m=S_m(Y,X)e_R^{\mathsf T}\). Reading the
\((R,R)\) entry of the inverse by cofactors gives
\begin{equation}
     \det J_m(X,Y)= S_m(Y,X)\det J_m^{\widehat R}(X,Y).
     \label{eq:rainbow-minor}
\end{equation}
This polynomial identity retains its meaning at a singular specialization.

\section{Polynomial certificates}
The supplement gives the two small polynomial equivalences at each of
orders four and five. At orders six, seven and eight the matrices have
dimensions 132, 429 and 1430. For each of these orders the certificate records
the sparse unitriangular multiplier; exact division reconstructs
the quotient. In the grading of Section~\ref{sec:groves}, a homogeneous
monomial of degree difference \(\delta\) specializing to \(x^j\)
must be \(X^{(j+\delta)/2}Y^{(j-\delta)/2}\).
The verifier checks both exponents for integrality and nonnegativity,
reconstructs each matrix entry from the two partition joins, and
checks the polynomial remainder of every entry. The degree argument
of Section~\ref{sec:degree} certifies the two determinants. At \(n=8\)
this verifies \(2{,}044{,}900\) entries in each ring. These finite
tests use no specialization at sampled integer values. A separate program evaluates
Corollary~\ref{cor:entries} and checks \(LG=DV\) coefficient by coefficient.
At order eight the quotient \(V\) has 675,204 nonzero entries taking 177 distinct
nonzero values over the two-variable ring for the supplied weights. The supplement also
checks the two shaded projector traces independently through six strands.
The polynomial certificates, verification scripts and both trace calculations
are available in the \href{https://github.com/Kimargin8/meander-ratio-and-smith/releases/tag/v2026.11.18}{open supplement}.

A separate check of \eqref{eq:rainbow-minor} starts from the two partition
joins and computes both determinants by elimination. For \(1\le m\le5\)
it verifies the polynomial identity over \(\mathbb Z[X,Y]\).
Indeed, every term on either side has the same difference of exponents
by the grading in Section~\ref{sec:groves}, since the black partition
of \(R\) has \(\lceil m/2\rceil\) blocks.
Setting \(Y=1\) therefore loses no coefficients. If \(C=C_m\) and
\(A=\sum_{\pi}|\pi|\), both sides have \(X\)-degree at most \(A\)
and are divisible by \(X^{C-1}\). Exact evaluations at respectively
3, 4, 8, 24 and 87 distinct positive values of \(X\) suffice.
The accompanying minor-check supplement also tests twelve finite-field
specializations at orders five through seven. At \(m=6\) over
\(\mathbb F_{101}\) the full determinant is 82 for both \((X,Y)=(2,3)\)
and \((3,2)\), while the respective principal minors are 35 and 2.

\paragraph{Declaration of generative AI and AI-assisted technologies
in the manuscript preparation process.}
The author used OpenAI Codex to develop proof arguments, search the
literature, write and test the verification programs, and draft and
revise the mathematical exposition.


\begin{thebibliography}{9}
\bibitem{ms} A. R. Miller and D. Stanton,
\emph{Orthogonal polynomials and Smith normal form},
Monatsh. Math. \textbf{187} (2018), 125--145.
\href{https://doi.org/10.1007/s00605-017-1082-6}{doi:10.1007/s00605-017-1082-6}.
In the accepted version see \S6.3.2, eqs.~(64)--(65), and
\S6.3.3, eq.~(69). The arXiv version states Conjecture~2 in \S6.2.3, eq.~(65).
\href{https://arxiv.org/abs/1704.03539}{arXiv:1704.03539}.
\bibitem{dahab} R. Dahab,
\emph{The Birkhoff--Lewis equations}, Ph.D. thesis,
University of Waterloo (1993), Theorems 1.8.1 and 2.5.2.
\bibitem{dfgg} P. Di Francesco, O. Golinelli and E. Guitter,
\emph{Meanders and the Temperley--Lieb algebra},
Commun. Math. Phys. \textbf{186} (1997), 1--59.
\href{https://arxiv.org/abs/hep-th/9602025}{arXiv:hep-th/9602025}.
\bibitem{webster} B. Webster, \emph{The existence of Jones--Wenzl projectors},
Appendix A to B. Elias and N. Libedinsky,
\emph{Indecomposable Soergel bimodules for universal Coxeter groups},
Trans. Amer. Math. Soc. \textbf{369} (2017), 3883--3910.
\href{https://doi.org/10.1090/tran/6754}{doi:10.1090/tran/6754}.
\bibitem{hazi} A. Hazi,
\emph{Existence and rotatability of the two-colored Jones--Wenzl projector},
Bull. Lond. Math. Soc. \textbf{56} (2024), 1095--1113.
\href{https://doi.org/10.1112/blms.12983}{doi:10.1112/blms.12983}.
\bibitem{chen} Q. Chen and J. H. Przytycki,
\emph{The Gram matrix of a Temperley--Lieb algebra is similar to the
matrix of chromatic joins}, \href{https://arxiv.org/abs/0806.0878}{arXiv:0806.0878}.
\bibitem{hlr} C. Howat, R. Laugwitz and M. Ray,
\emph{Temperley--Lieb categories with coloured regions and Jones--Wenzl projectors},
Bull. Lond. Math. Soc. \textbf{58} (2026), e70363.
\href{https://doi.org/10.1112/blms.70363}{doi:10.1112/blms.70363}.
\bibitem{morrison} S. Morrison, \emph{A formula for the Jones--Wenzl projections},
\href{https://arxiv.org/abs/1503.00384}{arXiv:1503.00384}.
\end{thebibliography}
\end{document}